\documentclass[graybox,reqn]{svmult}
\usepackage{mathptmx}       
\usepackage{helvet}         
\usepackage{courier}        
\usepackage{type1cm}        
\usepackage{makeidx}         
\usepackage{graphicx}        
\usepackage{multicol}        
\usepackage[bottom]{footmisc}
\usepackage{tikz} 

\usepackage{url}

\usepackage{amsfonts}
\usepackage{amsmath}
\usepackage{amssymb}

\renewcommand{\t}{\theta}
\renewcommand{\i}{\mathrm{i}}

\renewcommand{\div}{\mathrm{div}}
\newcommand{\ep}{\varepsilon}

\newcommand{\CC}{\mathbb{C}}

\newcommand{\RR}{\mathbb{R}}
\newcommand{\ZZ}{\mathbb{Z}}

\renewcommand{\E}{\mathcal{E}}
\renewcommand{\L}{\mathcal{L}}
\renewcommand{\H}{\mathcal{H}}

\newcommand{\V}{\mathcal{V}}

\newcommand{\Z}{\mathcal{Z}}

\begin{document}

\title*{Quantitative homogenisation for differential equations with highly anisotropic partially degenerating coefficients}
\titlerunning{Homogenisation for highly anisotropic fibres}
\author{Shane Cooper and Ilia Kamotski} 
\institute{Shane Cooper \at University College London, Gordon Street, London WC1E 6BT, UK. 
\and Ilia Kamotski\at University College London, Gordon Street, London WC1E 6BT, UK. \at \email{i.kamotski@ucl.ac.uk} }  

\maketitle

\abstract{
	We consider a non-uniformly elliptic second-order  differential operator with periodic coefficients that models composite media consisting of highly anisotropic cylindrical fibres periodically distributed in an isotropic background.  The degree of anisotropy is related to the period of the coefficients via  a `critical' high-contrast scaling.  In particular, ellipticity  is lost in certain directions as the period, $\ep$,  tends to zero.  Our primary interest is in the asymptotic behaviour of the resolvent of this operator in the limit of small $\ep$. 
	Two-scale resolvent convergence results were established for such operators in Cherednichenko, Smyshlyaev and Zhikov (Proceedings of The Royal Society of Edinburgh:Seciton A Mathematics. \textbf{136}(1), 87--114(2006)). In this work, we provide an asymptotic description of the resolvent and establish operator-type error estimates. 
Our approach adopts the general scheme of Cooper, Kamotski and Smyshlyaev (preprint available at https://arxiv.org/abs/2307.13151). However, we face new challenges such as a directional dependence on the loss of ellipticity in addition to a key `spectral gap' assumption of the above article only holding in a weaker sense. This results in an additional `interfacial' boundary layer analysis in the vicinity of each fibre to arrive at order-$\ep$ operator-type error estimates.
} 



\section{Introduction}
In this work, we analyse the asymptotic behaviour of the solution $u_\ep \in H^1(\RR^3)$ to the problem
	\begin{equation}\label{p.strong}
	\left.	\begin{aligned}
		&-\div\big( A_\ep(\tfrac{x}{\ep}) \nabla u_\ep(x) \big) + u_\ep(x) = f(x)  \ \text{for $f \in L^2(\RR^3)$,}\\
		&\ A_\ep(y) = \chi_1(y) I +  \chi_0(y) \mathrm{diag}\big(\ep^2,\ep^2,a(y_3)\big),  
				\end{aligned} \quad \right\}
		\end{equation}
		for small (period) parameter  $0<\ep<1$. Here,  $\square : = \big[-\frac{1}{2},\frac{1}{2}\big)^3$ is the periodic reference cell of $A_\ep$, $\chi_1=1-\chi_0$ and $\chi_0$ is the $\square$-periodic extension of the
	 characteristic function of the cylindrical set $F := B \times \RR$ where $B$ is the disk of radius $r<1/2$ centred at the origin. The $1$-periodic coefficient $a$ satisfies the usual well-posedness assumptions, i.e. there exists a positive constant $\nu>0$ such that
	 \[
	 v<a(s)<v^{-1}, \qquad \forall s \in [0,1).
	 \] 
	Such equations are one possible model for highly anisotropic periodic composite media with partial degeneracies, \cite{VPS}.

%
It was shown in \cite{ChSmZh} that the solution $u_\ep$ to \eqref{p.strong} two-scale converges, as $\ep \rightarrow 0$, to $u(x)+v(x,y')$ the solution of the two-scale homogenised limit problem
\begin{equation}\label{p4}
	\left. 	\begin{aligned}
		& - \div A^{\mathrm h} \nabla u - a^{\mathrm h} \partial_{x_3}^2 \big(  |B|u+\langle v\rangle \big)   +u+\langle v \rangle = f, &\qquad x\in \RR^3, \\
		&- \Delta_{y'} v - a^{\mathrm h}\partial_{x_3}^2  \left(u+ v \right) +u+v=  f, &\qquad  x\in \RR^3,  y' \in B, \\
	&\	v(x,y') = 0, &\quad x\in \RR^3,  y' \in \partial{B},
	\end{aligned} \quad \right\}
\end{equation}
where $\langle v \rangle(x) : = \int_B v(x,y') \, \mathrm{d}y'$ and $A^{\mathrm h}$ and $a^{\mathrm h}$ are given constant coefficients (defined below).
Note that the two-scale limit problem is defined on a subspace of $L^2(\RR^3\times \square)$ that is larger than the space on which the original problem is posed. As such, one has convergence in the two-scale topology rather than the standard topology, cf. the seminal works \cite{Ng},\cite{All} and \cite{Zh}. 

 In our recent work \cite{CoKaSm} we provide a general methodology to upgrade such two-scale convergence results to operator-type error estimates in the space on which the original problem is posed. 
 We begin by the standard procedure of quantitative homogenisation for periodic operators in the whole space
 to arrive at a family of equivalent problems on the unit cell $\square$. Namely, we sequentially apply, to $u_\ep$, the unitary rescaling $\Gamma_\ep : L^2(\RR^3) \rightarrow L^2(\RR^3)$, $\Gamma_\ep f(x) = \ep^{3/2} f(\ep x)$, then the unitary Gelfand transform $U: L^2(\RR^3) \rightarrow L^2( \square^* \times \square)$, for $\square^* =[-\pi,\pi)^3$ the dual cell and $U$ being the continuous extension of the mapping
\[
Uf(\t,y) = (2\pi)^{-3/2}\sum_{m\in \ZZ^3} f(y+m) e^{-\i \t \cdot  (y+m)}.
\] 
Then $u_{\ep,\t}(y): = \big(U \Gamma_\ep u\big) (\t,y)$ solves
\[
- \big( \div + \i \t \cdot \big) \ep^{-2}A_\ep (\nabla + \i \t) u_{\ep,\t} + u_{\ep,\t}= f_{\ep,\t} \quad y \in \square, \t \in \square^*, \quad 
\]
where $f_{\ep,\t}(y) : = \big(U \Gamma_\ep f\big) (\t,y) \in L^2(\square^*\times \square)$. Following \cite{CoKaSm}, we rewrite this problem as: for each $0<\ep<1$ and $\t\in \square^*$,
\begin{equation}\label{p0}
\left. \begin{aligned}
&\text{Find $u_{\ep,\t}\in H:= H^1_{per}(\square)$ the solution to}\\
&\ep^{-2}a_{\t}(u_{\ep,\t},\tilde{u}) + b_{\t}(u_{\ep,\t},\tilde{u} ) = c(f_{\ep,\t},\tilde{u} ), \quad \forall \tilde{u} \in H
\end{aligned}\quad \right\},
\end{equation}
for families of non-negative sesquilinear forms
 $a_\t, b_\t : H \times H \rightarrow \CC$ and $c: L^2(\square) \times L^2(\square) \rightarrow \CC$   given by the quadratic forms 
\begin{equation}\label{at}
	\begin{aligned}
	a_\t(u,u)=:	a_\t[u]  = 
		\int_{\square\backslash F} | (\nabla + \i \t) u(y)|^2\, \mathrm{d}y + \int_F a(y_3)|& (\partial_3 +\i \theta_3) u(y)|^2 \, \mathrm{d}y;
	\end{aligned}
\end{equation}
\begin{equation}\label{bt}
b_\t(u,u)=:	b_{\t}[u]  = \int_F  |(\nabla' +\i \theta') u|^2 + \int_\square |u|^2, \quad  \nabla'  = (\partial_1 , \partial_2 ), \ \t' = (\t_1,\t_2),
\end{equation}
and 
	\begin{equation*}
c[u]  = \int_\square  |u|^2. 
\end{equation*}
Notice that,  $a_\t + b_\t$ are a (uniformly equivalent) family of inner products on $H$ and so \eqref{p0} is well-posed.  

The goal now, is to provide leading-order asymptotics for solutions to the family of variational problems \eqref{p0}  that are uniform in $\t$, $\ep$ and the right-hand-side $f_{\ep,\t}$. Then one can invert the transforms $U$ and $\Gamma_\ep$ applied above to find uniform estimates for the problem \eqref{p.strong}. In \cite{CoKaSm}, such asymptotics were provided for variational families of the form \eqref{p0} with non-negative forms $a_\t$ that are (point-wise) coercive in $H$ (after projecting away from the kernel of $a_\t$). Amongst other things, this ensured that first order correction of $u_{\ep,\t}$ lived in $H$ which is crucial to the analysis in \cite{CoKaSm}. Another feature of \cite{CoKaSm} is that the coercivity constant degenerates quadratically in $\t$. 

Here, the form   $a_\t$ (cf. \eqref{at}) is clearly not coercive (up to the kernel) in $H^1_{per}(\square)$. Indeed,  consider
\begin{equation}\label{V}
V  
 = 	\{ v \in H^1_{per}(\square) : \text{$v$ constant in $\square \backslash F$ and  $v$ independent of $y_3$ in $F$} \}.
\end{equation}
Clearly, $a_\t$ is zero on $V$  when $\t_3=0$. Now, denote  the closed subset $W $ of $H^1_{per}(\square)$ that is the orthogonal complement of $V$  with respect to the  inner product $a_0 + b_0$ (which is an equivalent $H^1$ inner product).  Then, it is easy to see that $a_\t$ is positive on $W$. The key assumption of \cite{CoKaSm}	 (up to some details) was that $a_\t$ is bounded from below by $b_\t$ on $W$. This `spectral-gap' assumption, which holds for many examples of interest, cf. \cite{IKVPS} and \cite{CoKaSm}, fails in the present example (i.e. for forms \eqref{at} and \eqref{bt}). However, one has the following weaker coercivity estimate:  there exists a positive constant $\gamma$ such that
	\begin{align*}
	a_\t[w] &\ge  \gamma c[w], \quad \forall w \in W, \ \forall \t \in \square^*. 
\end{align*}
Overcoming this absence of full coercivity is one of the main challenges of this current work as it leads to  new analysis not present in \cite{CoKaSm}, in particular, we need to investigate internal boundary-layer effects near each fibre.

Another difficultly is related to the directional dependence of the elliptic degeneracy, which is a feature in various highly anisotropic high-contrast composite media (see for example \cite{Co}), and the failure of the quadratic degeneracy assumption of \cite{CoKaSm}. Loosely speaking, in the present context, this assumption would be satisfied if one could demonstrate that $a_\t[u]$ is bounded from below by $|\t|^2 c[u]$ for all $u$. This is clearly not the case. However, for \eqref{at}, one has the following `directionally-dependent' quadratic condition: there exists a positive constant $\gamma^*$ such that
 \begin{equation*}
 	a_\t[u] \ge \gamma^* \Big( |\t_3|^2 c[\chi_0 u] + |\t|^2 c[ \chi_1 u] \Big), \quad \forall u \in H^1_{per}(\square), \forall \t \in \square^*. 
 \end{equation*}
 Such a weaking the quadratic gap assumption of \cite{CoKaSm} present significant challenges. These challenges are too significant to achieve results on the same level of  generality as in our previous work \cite{CoKaSm}.  However, we were able to succeed for problem \eqref{p.strong} by exploiting problem-specific features.

\section{Results}
Our first result concerns uniform asymptotics of \eqref{p0} with operator-norm error estimates. 
Upon representing $V$, cf. \eqref{V},  as the direct sum $V = Z_0 \dot{+} Z_1$, where
\[
\begin{aligned}
Z_0 :=\CC,\ Z_1 : = \{ v \in H^1_{per}(\square) \, : \,  \text{$v=0$  in $\square \backslash F$ and  $v$ independent of $y_3$ in $F$}   \}, 
\end{aligned}
\]
we introduce the following problem on $V$:
\begin{equation*}
	\left. \begin{aligned}
		&\text{Find $(z_0,z_1) \in Z_0 \times Z_1$ the solution to}\\
		&\ep^{-2} A^{\mathrm h} \t \cdot \t\, z_0 \overline{\tilde{z}_0}  + a^{\mathrm h}\ep^{-2}|\t_3|^2   \big( z_0+ z_1,\tilde{z}_0+ \tilde{z}_1\big)_{L^2(B)} \\
		& \hspace{.17\linewidth}+  b_{0}(z_0+ z_1,\tilde{z}_0+ \tilde{z}_1 ) 
		 = c( \E_\t f_{\ep,\t},\tilde{z}_0 + \tilde{z}_1 ) \quad \forall \tilde{z}_0 \in Z_0,  \tilde{z}_1 \in Z_1
	\end{aligned} \quad \right\}.
\end{equation*}
Here, $\E_\t : L^2(\square) \rightarrow L^2(\square)$ is the unitary operator given by multiplication by the function $\exp(\i \t' \cdot y')$, $a^{\mathrm h}: = \left(\int_0^1 a^{-1} \right)^{-1}$ is the one-dimensional homogenised matrix and $A^{\mathrm h}$ is the constant, symmetric $3\times 3$  perforated-fibre homogenised matrix given, in \cite{ChSmZh}, as 
\[
A^{\mathrm h} : = |\square\backslash F| e_3\otimes e_3 +  \sum_{\alpha,\beta=1}^2 e_\alpha \otimes e_\beta\int_{\square'\backslash B} (\partial_{\beta} N_\alpha + \delta_{\alpha \beta}),
\] 
for  $N_\alpha : H^1_{per}(\square'\backslash B) \rightarrow \RR$, $\square' :=\big[-\frac{1}{2},\frac{1}{2}\big)^2$, $\alpha=1,2$,  that solve the `cell problem'
\begin{equation*}
\left. \begin{aligned}
&\int_{\square'\backslash B} (\nabla' N_\alpha + e_\alpha ) \nabla' \phi = 0, \quad \forall \phi \in H^1_{per}(\square'\backslash B)\\
& \int_{\square' \backslash B} N_\alpha= 0
\end{aligned} \quad \right\}.
\end{equation*}
Here $e_k$, $k=1,2,3$, are the canonical Euclidean basis vectors.


Our first main result is the following uniform estimate: 
\[
\| u_{\ep,\t} - \E_\t^*(z_0 +  z_1) \|_{L^2(\square)} \lesssim \ep \| f_{\ep,\t}\|_{L^2(\square)}, 
\]
which can be rewritten in operator-theoretic language as follows:
\begin{theorem}\label{thm1.2}
	There exists a positive constant $C$ such that the inequality
	\[
	\| (\mathcal{L}_{\ep,\t}+I)^{-1} - \E_\t ^*(\mathbb{L}_{\t/\ep}+I)^{-1} P_0 \E_\t \|_{L^2(\square) \rightarrow L^2(\square)} \le C \ep
	\]
holds for all  $0<\ep<1$ and $\t\in \square^*$.  
\end{theorem}
The operators in the above theorem are the following: for each $\ep$ and $\t$, $(\mathcal{L}_{\ep,\t}+I)$ is the unbounded self-adjoint operator in $L^2(\square)$ generated by problem \eqref{p0}; now
denote $\H_0$ to be the closure of $V$ in $L^2(\square)$ and $P_0 : L^2(\square) \rightarrow \H_0$ the orthogonal projection; then, finally, for each $\xi \in \RR^3$,  $(\mathbb{L}_\xi+I)$ is the unbounded self-adjoint operator in $\H_0$ generated by the inner product
\[
\mathbb{S}_\xi(z_0+z_1,\tilde{z}_0+\tilde{z}_1) : =  A^{\mathrm h} \xi \cdot \xi z_0 \overline{\tilde{z}_0}  + a^{\mathrm h}|\xi_3|^2  \big(z_0+  z_1,\tilde{z}_0+ \tilde{z}_1\big)_{L^2(B)} +  b_{0}(z_0+ z_1,\tilde{z}_0+ \tilde{z_1} ).
\] 

Each of the non-negative operators  $\mathcal{L}_{\ep,\t}$ and $\mathbb{L}_{\xi}$ have compact resolvents. Consequently, Theorem \ref{thm1.2} readily implies uniform approximations for the associated eigenvalues. Indeed, for $\{\lambda_{\ep}^{(k)}(\t) \}_{k=1}^\infty$ and $\{\Lambda^{(k)}(\xi) \}_{k=1}^\infty$ denoting the eigenvalues of $\mathcal{L}_{\ep,\t}$ and $\mathbb{L}_\xi$ respectively, labelled in increasing order and repeated according to multiplicity, one has the following result: 
\begin{theorem}\label{thm1.3}
	For each $k \in \mathbb{N}$, there exists a constant $C_k>0$ such that
\[
\left| \lambda_{\ep}^{(k)}(\t) - \Lambda^{(k)}\left(\tfrac{\t}{\ep}\right) \right| \le C_k \ep, \quad   \forall \t \in \square^*, \, \forall\, 0<\ep<1.
\]
\end{theorem} 
 Theorem \ref{thm1.2} enables one to approximate the solution of the original problem \eqref{p.strong} in terms of the two-scale homogenized limit problem \eqref{p4}. To this end, let us define the Hilbert space 
\[
\H : = L^2(\RR^3;\mathcal{H}_0),
\]
$\mathcal{P} : L^2(\RR^3\times \square) \rightarrow \H$ the orthogonal projection,  and define the dense subspace
\[
\V : = H^1(\RR^3) \dot{+} \Z_1, \quad  \Z_1 : = \{ v \in  L^2(\RR^3,Z_1) \, : \, \partial_{x_3} v \in L^2(\RR^3\times \square)\},
\]
of $\H$. Finally, let $\L_0+I$ be the unbounded self-adjoint operator defined by the closed quadratic form
\begin{equation*}
	 	\begin{aligned}
& Q[u+v] = \int_{\RR^3} \Big(	 A^{\mathrm h} \nabla u \cdot \overline{\nabla u}  + \  a^{\mathrm h}  \|\partial_{x_3}( u+  v)\|^2_{L^2(B)}   +  b_{0}[u +v] \Big)\, \mathrm{d}x,
	\end{aligned} 
\end{equation*}
with form domain $\V$.

To provide error estimates in $\mathcal{L}(L^2(\RR^3))$ between the resolvents of  $\mathcal{L}_\ep$ and $\mathcal{L}_0$,  one needs to address the fact that the operators are defined on different function spaces. We  overcome this issue via a `two-scale interpolation' operator $J_\ep : L^2(\RR^3) \rightarrow L^2(\RR^3 \times \square)$, introduced in \cite{CoKaSm} (see \cite{Well} for a similar object)  given by 
\[
J_\ep : =   \Gamma_\ep^{-1}\mathcal{F}^{-1}\chi \E_\t^*U\Gamma_\ep.
\] 
Here, $\chi : L^2(\square^* \times \square) \rightarrow L^2(\RR^3 \times \square)$ extends functions by zero outside of $\square^*$; $\mathcal{F}$ is the Fourier transform acting on the first variable (and $\Gamma_\ep$, $\Gamma_\ep^{-1}$ are understood to act on the first variable).  It was established in \cite{CoKaSm} that $J_\ep$  is almost-unitary in the following sense:
 \begin{equation*}
 	J_\ep^* J_\ep = I, \quad \text{and}\quad J_\ep J^*_\ep  \rightarrow  I \ \text{strongly as $\ep \rightarrow 0$.} 
 \end{equation*}
Then, the following result holds: 
  \begin{theorem}\label{thm2}
 	There exists a constant $K>0$ such that the  inequality
 	\[
 	\| (\mathcal{L}_\ep+I)^{-1} - J_\ep^* (\mathcal{L}_0+I)^{-1} \mathcal{P} J_\ep \|_{\mathcal{L}(L^2(\RR^3))} \le K \ep,
 	\]
 	holds for all $0<\ep<1$. 
 \end{theorem}
 Our results imply comparisons between the spectral properties of  $\mathcal{L}_\ep$ and $\mathcal{L}_0$.  For example, upon noting the fibre decomposition $\big( \mathcal{F}\mathcal{L}_0\mathcal{F}^{-1} \big)(\xi) =  \mathbb{L}_\xi$,  one sees that Theorem \ref{thm1.3} provides uniform estimates between the spectral band functions of the operators $\mathcal{L}_\ep$ and $\mathcal{L}_0$.
 	
The proof of the announced results will be appear elsewhere. 

\begin{acknowledgement}
The authors would like to thank Professor Valery Smyshlyaev for bringing this problem to our attention. 
\end{acknowledgement}

\end{document}